\documentclass[12pt]{article}
\usepackage{amsmath}
\usepackage{amssymb}
\usepackage{amsthm}
\usepackage{amscd}
\usepackage{amsfonts}
\usepackage{graphicx}%
\usepackage{fancyhdr}

\graphicspath{ {Images/} }

\usepackage{tikz}
\usepackage{verbatim}

\usetikzlibrary{matrix}
\usepackage{color, colortbl}

\begin{document}

\title{Math Course Redesign in a Private Four Year Hispanic Serving Institute to Address DEI Issues}

\author{C. Chang, Z. Chen}

\date{}

\maketitle

\section{Background}
Our institution is a private four-year college in the New York metro area, with total student enrollment of 10,577 in 2019-2020. The College is a Hispanic Serving Institute with 70\% female students. The enrollment of calculus classes has been consistently low despite a relatively large population of STEM majors, moreover, the female student’s ratio in those classes is significantly less than 10\%. Compounded by the fact that commuter students in the College made up 60\% of the population, students in general do not have a strong sense of community. This lack of inclusiveness is more significant for the limited number of women students taking calculus classes.
   
We recently organized a math symposium\footnote{https://www.mercy.edu/liberal-arts/math-symposium} on redesigning math for student engagement and success, to share experiences among some regional institutions of similar type. In our own case, we identified three most challenging points related to diverse, equitable, and inclusive (DEI) issues. First, the majority of our students entering the College lack the math skills essential to success in Calculus, as basic as College Algebra, some others have a multi-year gap after graduating high school. Almost all but a few STEM students must start from College Algebra before they can move on to Precalculus and then Calculus. Secondly, we noted that many students who planned to pursue STEM dropped out of their majors because they couldn’t obtain the required grade in College Algebra to move forward. This is one of the main reasons that the enrollment of calculus classes is consistently low. Lastly, a large portion of basic math classes are taught by adjunct instructors, the turnover ratio among adjunct instructors is not small. One such consequence is that many students don’t have equitable learning experiences and some students are still struggling with College Algebra even in the calculus class.

In order to have a robust and successful calculus sequence, we need to redesign our basic math curriculum and build a strong foundation for students, so that they can channel through to the calculus class. Especially for our underrepresented population like first-generation college students or female students in the STEM fields, they may not have received the attention, support, and resources in their high school education, so the academic preparedness gap for them are detrimental in college math courses and could well be the reason for eventual college drop off.

\section{Smart Math Initiative}
The name of our program to solve DEI issues is called Smart Math Initiative and its benefit is in twofold. For one, it is designed to provide a standardized engaging teaching pedagogy across all basic math classes, diminish the vast differences among adjunct instructors, guarantee the rigor and consistency of the curriculum, and accommodate each student with an individualized and self-paced study plan. For another, the Smart Math Initiative provides students, especially underrepresented students like women, an abundance of support and encouragement inside and outside of classrooms. There are three main components.

\subsection{Embedded Peer Tutor (EPT)}

The Embedded Peer Tutors (EPTs) are in-class student teaching assistants who have previously taken and received an A in a math class equivalent to or higher than College Algebra. They are currently enrolled undergraduate students of any major, as early as second semester in the freshman year. Unlike working in tutoring centers, EPTs work directly with students inside classrooms, in conjunction with and complement to the instructors.

We usually hire 10 to 20 EPTs each semester, depending on their schedules and availability. Each semester, all the EPTs must go through a training before they work inside the classrooms. By the nature of the emporium model which will be elaborated in the next component, a class is usually divided up into small groups or even individuals and have separate discussions or learning activities from time to time. Therefore, EPTs are supervised and directed by instructors to assist students, in the form of mini lecturing a learning topic, facilitating a group discussion, or helping students self-reflect learning strategies and form a good working habit.
 
We found, through the end of semester student evaluations, that the emotional support and encouragement provided by EPTs are indispensable to peer students, especially the struggling students who have a lot of fear in learning mathematics. The role of EPTs in creating an inclusive learning environment is irreplaceable by instructors.

An interesting observation from the past, the majority of our EPTs are female, most of them are STEM major students. Many EPTs eventually enrolled into a 5-year master program in teaching education degrees in STEM.

\subsection{Controlled Emporium Model in Smart Math Lab}

All the Smart Math Initiative classes take place in the designated computer labs, called Smart Math Labs, to create an equitable and regulated classroom environment. These labs have open floor concepts, with a maximum capacity of 30 students. The design of the labs is structured to welcome class discussions, nurture a support environment, and create bonds among students.
 
In terms of the course contents, they are divided up into modules, delivered through a learning platform like MyOpenMath or ALEKS. Outside the classrooms, students progress through the program at their own pace, by a standardized benchmark timeline. With self-paced learning, students usually dive into different topics without realizing the connection and the big picture of them, and this will be addressed inside the classrooms.
 
In a typical 2-hour class meeting, an instructor will motivate and teach between half an hour and an hour, on a specific topic outlined in the benchmark schedule. The lecture will provide additional examples as well as connections between different topics. The rest of the class meeting will plan out as individual interactions or small group discussions between the instructor and the students, assisted by EPTs. 

Following closely the Instructional Practices Guide\footnote{Abell, M., Braddy, L., Ensley, D., Ludwig, L., \& Soto-Johnson, H. (2017). MAA Instructional Practices Guide.}, the department regularly organizes workshops on best teaching strategies. Some high standard and consistent teaching materials including well-designed group discussion worksheets are developed and updated for all instructors. These efforts are made to guarantee that students receive maximum support, and instructors, full-time or part-time, have the same standard in delivering the courses.

\subsection{Coordination and Technology}

The College designated a full-time faculty within the math department as the Initiative Coordinator, the key responsibility is coordinating with different departments. We meet with the institutional research for data and performance analysis, we work with the center for teaching and learning on workshops and class observations to support adjunct instructors, and our department also works with the tutoring center to recruit and supervise EPTs.  

The textbooks, videos, and other readings for the courses are all available electronically online, conveniently accessible through mobile devices as well as computers. Instructors receive electronic reports on students’ weekly performances, the coordinator monitors the learning statistics across sections. A typical weekly report includes the topics that the majority of the students are struggling with, the total number of hours that students spent on the learning platform, and the average grade of the weekly assignment. Through reading the reports prior to class meetings, instructors can identify the weaker students and help them more efficiently during class discussions. 

Finally, the coordinator maintains an online forum, to create a virtual community among instructors. A shared site for instructors and EPTs is regularly updated with teaching notes, pedagogical articles, extra exercises and detailed solutions, etc.

\section{Assessment}
One of the key metrics in measuring the success of the Smart Math Initiative is by the completion rate. The two courses redesigned are Math115 (Math for Liberal Arts) and Math116 (College Algebra). Our focus in this article is on Math116, as it is required by STEM majors and its population is the potential pool for future calculus students.

\begin{figure}[h]
    \centering
    \includegraphics[scale=0.6]{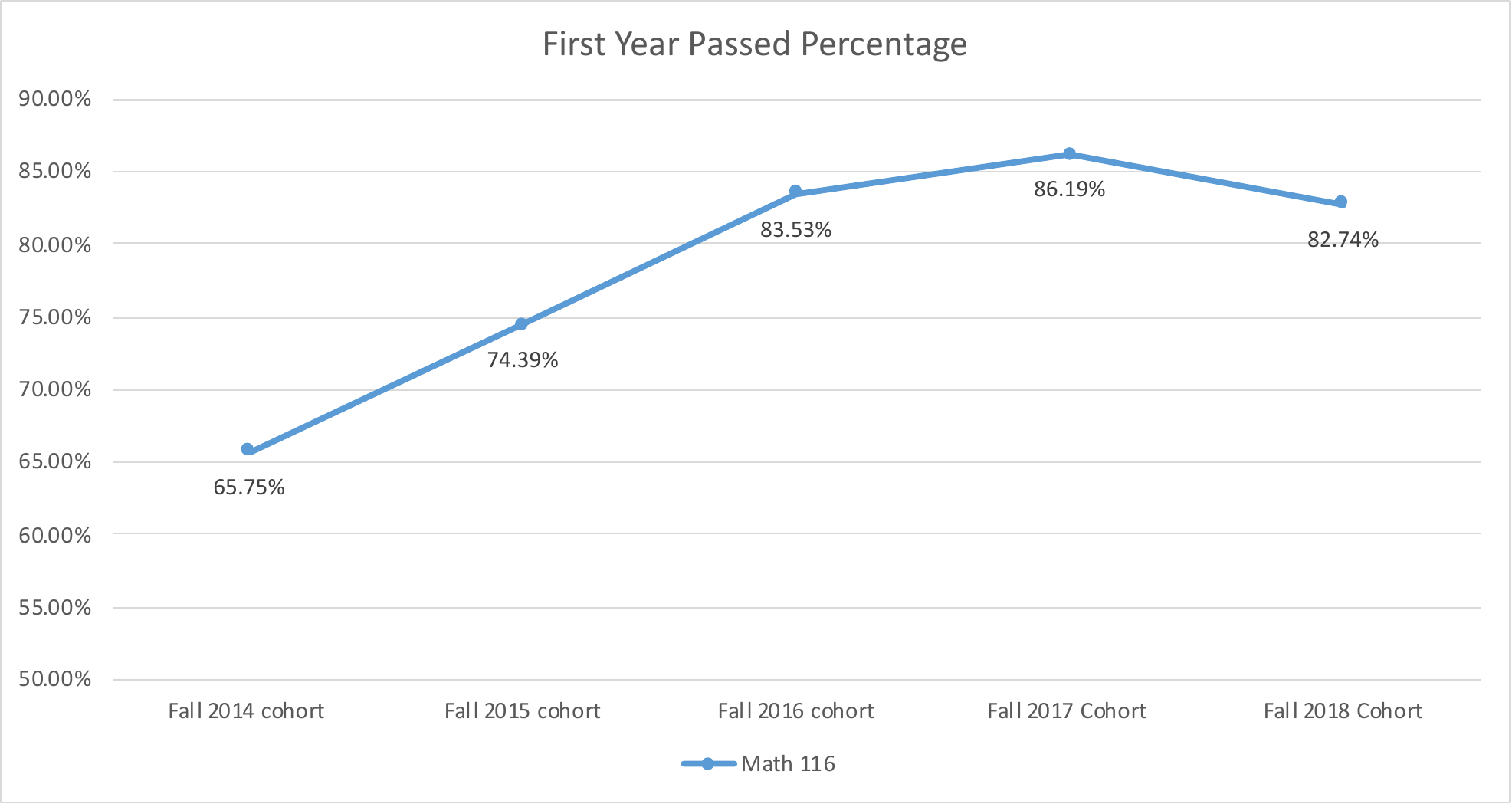}
    \caption{Math116 First Year Passed Percentage}
    \label{F:P}
\end{figure}

Table \ref{T:P} lists the completion rates for five different cohorts. By a cohort, we mean all freshman students enrolled in a fall semester, excluding continuing and transfer students. We track each cohort up to three years, since some students decided not to take Math116 in the first year or second year, some of them waited until the third year or beyond. The passed percentage is based on the number of students who took and passed Math116 with a grade of C or better, out of the total number of students in a cohort. Figure \ref{F:P} shows the trending of the first year passed percentage. The Initiative, which began in Spring 2016, is the driving force of the steady increase.

\begin{table}
    \centering
   \resizebox{\textwidth}{!}{
         
    \begin{tabular}{l|r r r r r}
    & F14 cohort &F15 cohort&F16 cohort&F17 cohort&F18 cohort \\
\hline
Total Students     & 619 &785 & 662&659&643\\
Passed within first year     & 407 &584 & 553&568&532\\
\rowcolor{yellow}
Passed Percentage     & 65.75\% &74.39\% & 83.53\%&86.19\%&82.74\%\\
Passed within two years     & 115 &116 & 26&15&\\\rowcolor{yellow}
Passed Percentage ( yr 1,2)     & 84.33\% &89.17\% & 87.46\%&88.47\%&\cellcolor{white}\\
Passed within 3 years     & 33 &10 & 1&&\\
\rowcolor{yellow}
Passed Percentage ( yr 1,2,3)     & 89.66\% &90.45\% & 87.61\%&\cellcolor{white}&\cellcolor{white}
\end{tabular}}

\caption{ Three-year Math116 completion rates comparison: F14 and F15 cohorts before Smart Math; F16, F17, and F18 cohorts after Smart Math}

    \label{T:P}
\end{table}

Next, we want to examine how the EPTs help the classes in passing rates. We have a large percentage of female EPTs, in Spring 2019, we hired 13 EPTs, out of which 10 are female. The female EPTs ratios in Fall 2019 and Spring 2020 are 8 out of 14 and 6 out 8, respectively. Table \ref{T:2} compares the pass rates of Math116 in Spring 2019, categorized by sections with EPTs and without. There is a clear contrast in the DFW rates, and the sections with EPTs significantly outperformed those without EPTs. 

\begin{table}
 \centering

    \begin{tabular}{l|r|r c | l|r|r}
     \cline{1-3}
     \cline{5-7}
\rowcolor{lightgray}      Numbers of students   &  A-C & DFW&\multicolumn{1}{l}{\cellcolor{white}}  &Percentages& A-C& DFW  \\ \cline{1-3}
     \cline{5-7}
    Not-EPTs&36&58&\multicolumn{1}{l}{}&Not-EPTs&38.30\%&61.70\%  \\\cline{1-3}
     \cline{5-7} 
    EPTs&198&118&\multicolumn{1}{l}{}&EPTs&62.66\%&37.34\%\\\cline{1-3}
     \cline{5-7}
    \end{tabular}
    \caption{Spring 2019 Math116 EPTs Analysis}
    \label{T:2}
\end{table}

\begin{figure}[h]
    \centering
    \includegraphics[scale=0.6]{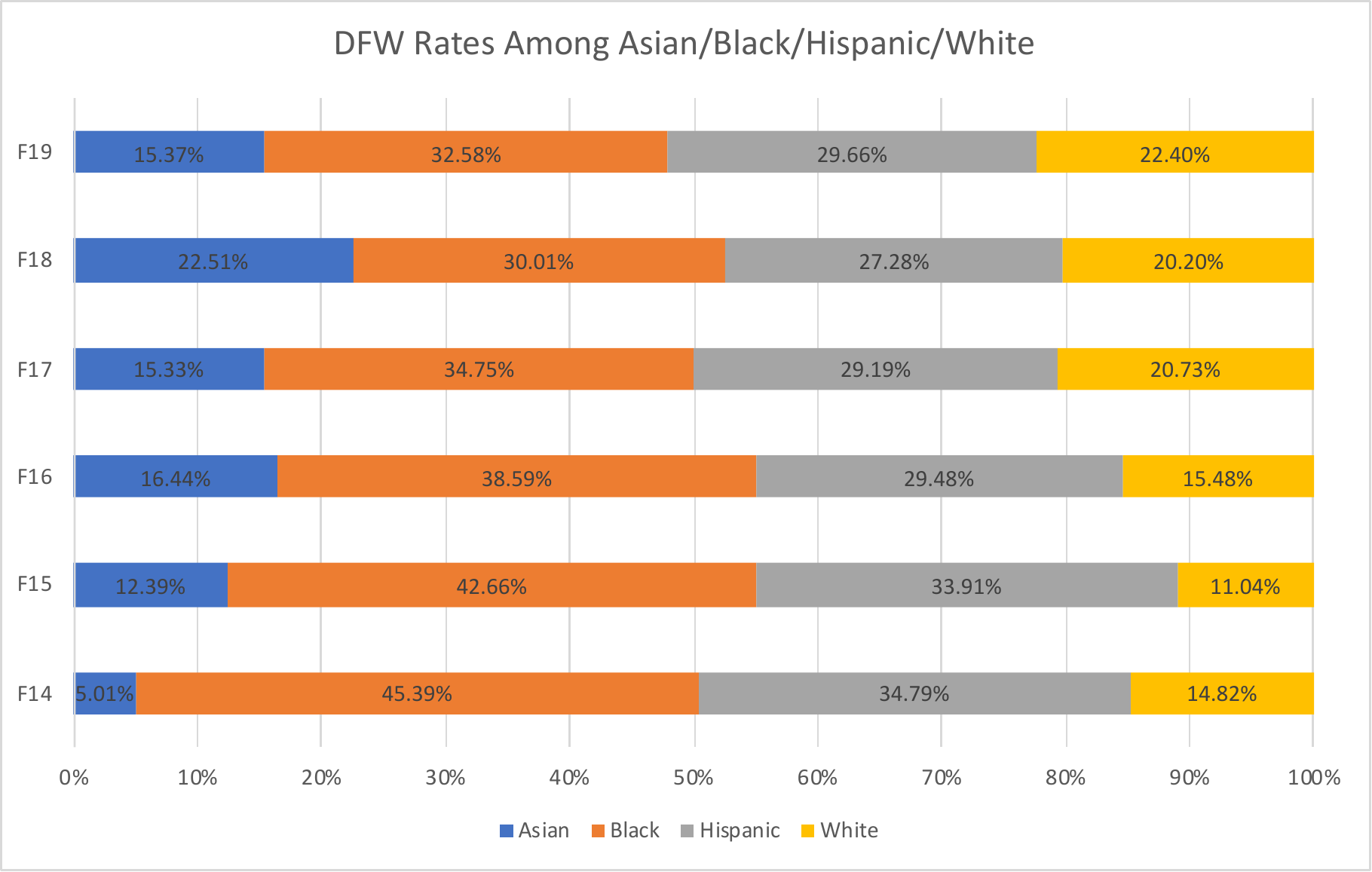}
    \caption{Math116 DFW rates among Asian/Black/Hispanic/White}
    \label{F:E}
\end{figure}

Regarding the performance of different ethnicity groups, Table \ref{T:E} provides the consolidated grade distributions within each group for the last six consecutive fall semesters. To make a direct comparison, we normalized the number of students in each group to 1 and focused on the DFW students. In the absolute ideal scenario, where all ethnicity groups perform equally well, the percentage of each group among the DFW students will be 25\%. Figure \ref{F:E} displays the distributions, and it clearly indicates that since the Initiative fully implemented in Spring 2016, the disparity among different ethnicity groups are reducing. 

\begin{table}
    \centering
    \resizebox{\textwidth}{!}{
    \begin{tabular}{l|r r r r r c l|r r r r r}
     \cline{1-6}
\cline{8-13}
\rowcolor{lightgray}
     Fall 14 &Asian& Black&Hispanic&White&Overall&\multicolumn{1}{l}{\cellcolor{white}}&Fall 15& Asian& Black&Hispanic&White&Overall\\
     \cline{1-6}
\cline{8-13}
     A-C&96.77\%&70.75\%&77.58\%&90.45\%&80.68\%&&A-C&90.32\%&66.67\%&73.50\%&91.37\%&78.84\%\\
\cline{1-6}
\cline{8-13}
\rowcolor{yellow}
DFW& 3.23\%&29.25\%&22.42\%&9.55\%&19.32\%&\cellcolor{white}&DFW&9.68\%&33.33\%&26.50\%&8.63\%&21.16\% \\
\cline{1-6}
\cline{8-13}
\multicolumn{1}{l}{}

\\
 \cline{1-6}
\cline{8-13}
\rowcolor{lightgray}
     Fall 16 &Asian& Black&Hispanic&White&Overall&\multicolumn{1}{l}{\cellcolor{white}}&Fall 17& Asian& Black&Hispanic&White&Overall\\
     \cline{1-6}
\cline{8-13}
     A-C&82.76\%&59.54\%&69.09\%&83.77\%&70.41\%&&A-C&87.50\%&71.67\%&76.20\%&83.10\%&75.83\%\\
\cline{1-6}
\cline{8-13}
\rowcolor{yellow}
DFW& 17.24\%&40.46\%&30.91\%&16.23\%&29.59\%&\cellcolor{white}&DFW&12.50\%&28.33\%&23.80\%&16.90\%&24.17\% \\
\cline{1-6}
\cline{8-13} 
\multicolumn{1}{l}{}
\\

 \cline{1-6}
\cline{8-13}
\rowcolor{lightgray}
     Fall 18 &Asian& Black&Hispanic&White&Overall&\multicolumn{1}{l}{\cellcolor{white}}&Fall 19& Asian& Black&Hispanic&White&Overall\\
     \cline{1-6}
\cline{8-13}
     A-C&75.00\%&66.67\%&69.70\%&77.56\%&70.46\%&&A-C&83.33\%&64.68\%&67.84\%&75.71\%&69.41\%\\
\cline{1-6}
\cline{8-13}
\rowcolor{yellow}
DFW& 25.00\%&33.33\%&30.30\%&22.44\%&29.54\%&\cellcolor{white}&DFW&16.67\%&35.32\%&32.16\%&24.29\%&30.59\% \\
\cline{1-6}
\cline{8-13}

\end{tabular}
}
\caption{Math116 consolidated grade distribution within  ethnicity groups} 
\label{T:E}
\end{table}

\begin{figure}
    \centering
    \includegraphics[scale=.55]{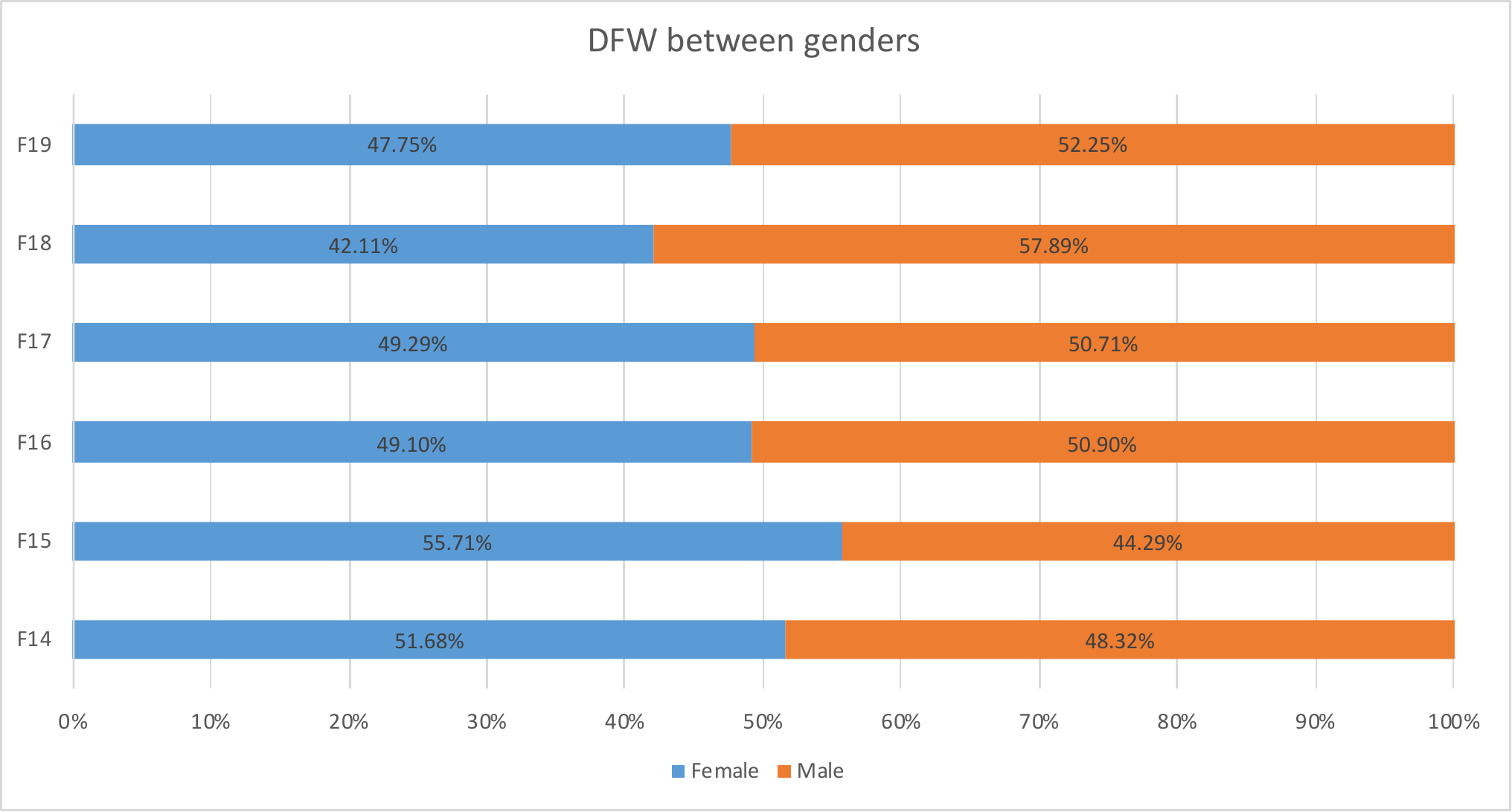}
    \caption{Math116 DFW rates between genders}
    \label{F:G}
\end{figure}

Similarly, Table \ref{T:G} lists the consolidated grade distributions within each gender over the past six fall semesters, and Figure \ref{F:G} displays the direct comparison between genders. For the reference, in an ideal situation where male and female students perform equally well, the percentage of male or female among DFW students should be 50\%. We observe in Figure \ref{F:G} that the percentage of females among DFW students is gradually shrinking below 50\%.

\begin{table}[h]
    \centering
    \scalebox{0.75}{
    \begin{tabular}{l|r r  c l|r r }
     \cline{1-3}
\cline{5-7}
\rowcolor{lightgray}
     Fall 14 &Female& Male&\multicolumn{1}{l}{\cellcolor{white}}&Fall 15& Female& Male\\
     \cline{1-3}
\cline{5-7}
     A-C&76.52\%&78.05\%&&A-C&75.97\%&80.90\%\\
\cline{1-3}
\cline{5-7}
\rowcolor{yellow}
DFW& 23.48\%&21.95\%&\cellcolor{white}&DFW&24.03\%&19.10\% \\
\cline{1-3}
\cline{5-7}
\multicolumn{1}{l}{}

\\
 \cline{1-3}
\cline{5-7}
\rowcolor{lightgray}
     Fall 16 &Female& Male&\multicolumn{1}{l}{\cellcolor{white}}&Fall 17& Female& Male\\
     \cline{1-3}
\cline{5-7}
     A-C&72.14\%&71.12\%&&A-C&23.65\%&24.34\%\\
\cline{1-3}
\cline{5-7}
\rowcolor{yellow}
DFW& 27.86\%&28.88\%&\cellcolor{white}&DFW&23.65\%&24.34\%\\
\cline{1-3}
\cline{5-7} 
\multicolumn{1}{l}{}
\\

 \cline{1-3}
\cline{5-7}
\rowcolor{lightgray}
     Fall 18 &Female& Male&\multicolumn{1}{l}{\cellcolor{white}}&Fall 19& Female& Male\\
     \cline{1-3}
\cline{5-7}
     A-C&74.16\%&64.48\%&&A-C&70.80\%&68.04\%\\
\cline{1-3}
\cline{5-7}
\rowcolor{yellow}
DFW& 25.84\%&35.52\%&\cellcolor{white}&DFW&29.20\%&31.96\%\\
\cline{1-3}
\cline{5-7}

\end{tabular}
}

 \caption{Math116 consolidated grade distribution within  genders}
    \label{T:G}
\end{table}

\section{Implications and Recommendations}

Since we started the Smart Math Initiative, the EPTs program has always been a core component. For students, they not only look up to the EPTs as role models, but also, they see in themselves that they can be as successful as the EPTs. Another unexpected outcome in EPTs is that a large percentage of them are female students, around 75\% in two of the most recent three semesters. 

From a financial point of view, EPTs receive stipends for their work and relieve their financial pressure at some level. The close involvement of EPTs in the classrooms, as well as their achievements in helping peer students have strengthened the inclusiveness and community bonds among students, and hence improve the retention rate.
  
The emporium model has been explored for many years across the states, we make some small tweaks in implementing to fit our unique needs. For instance, a not small portion of our student bodies are nontraditional college students, some have full-time jobs outside the school while maintaining a student status. We accommodate these students by not entirely flipping the classroom, instead, we use a combination of engaging lecturing and self-paced active learning. In order for our controlled emporium model to be successful, we spent a lot of efforts in workshops and training for instructors, especially part-time instructors.

Despite the initial success of the Smart Math Initiative, there are areas we are actively seeking to improve in the future. There are ideas about tackling these issues and the prospect of the Smart Math Initiative has become more exciting than it has ever been before.

\vspace{2cm}

\textbf{C. Chang}, Assistant Professor, Department of Math \& CS, Mercy College.

\textbf{Z. Chen}, Professor, Department of Math \& CS, Mercy College.

\end{document}